\input amstex
\magnification=1200
\documentstyle{amsppt}
\loadbold
\redefine\!{\kern-.075em}

\predefine\accute{\'}
\redefine\'{\kern.05em{}}
\def\bbigskip{\vskip 24pt plus 8pt minus 8pt}
\def\qed{\hfill$\ssize\square$}
\def\<{\langle}
\def\>{\rangle}
\def\leqs{\leqslant}
\def\geqs{\geqslant}
\def\Tr{\operatorname{Tr}}

\def\divergence{\operatorname{div}}
\def\supp{\operatorname{supp}}
\def\const{\operatorname{const.}}
\def\dotprod{\'\raise.25ex\hbox{$\sssize\bullet$}\'}
\def\F{\Cal F}
\def\H{\Cal H}
\def\J{\Cal J}
\def\a{\alpha}
\def\b{\beta}
\def\d{\delta}
\def\s{\sigma}
\def\v{\varphi}
\def\btheta{{\boldsymbol\theta}}
\def\bphi{{\boldsymbol\phi}}
\def\nab#1{\nabla\kern-.2em\lower.8ex\hbox{$\ssize#1$}\'}
\PSAMSFonts

\NoBlackBoxes
\nologo
\TagsOnRight
\rightheadtext{Harmonic Unit Vector Fields in 3-Space}
\pagewidth{16.2truecm}
\pageheight{24.3truecm}
\voffset=-0.5truecm
\linespacing{1.5}
\topmatter
\title 
Bending and Stretching Unit Vector Fields in Euclidean and Hyperbolic
3-Space
\endtitle
\author 
C\. M\. Wood
\endauthor
\address 
Department of Mathematics, University of York, Heslington, York
Y010 5DD, UK.
\endaddress
\email 
cmw4\@york.ac.uk 
\endemail
\dedicatory
For Esther.
\enddedicatory
\abstract
New examples of harmonic unit vector fields on hyperbolic 3-space are
constructed by exploiting the reduction of symmetry arising from the
foliation by horospheres.  This is compared and contrasted with the
analogous construction in Euclidean 3-space, using a foliation by planes, 
which produces some new examples of harmonic maps from 3-dimensional
Euclidean domains to the 2-sphere.  Finally, the harmonic unit vector field tangent to a parallel family of hyperbolic geodesics is shown to be unstable, by
constructing a class of compactly supported energy-decreasing variations.
All examples considered have infinite total bending.
\endabstract
\keywords
Harmonic section, harmonic map, harmonic unit vector field, harmonic 
function, sine-Gordon equation, pendulum equation, horospheres, energy,
bending, unstable, stretchable
\endkeywords
\subjclassyear{2000}
\subjclass
53C43, 58E20
\endsubjclass
\endtopmatter
\document
\head
1. Introduction
\endhead
A smooth unit vector field $\s$ on an orientable Riemannian manifold $(M,g)$ is said to be {\sl harmonic\/} if $\s$ is a harmonic section of the unit tangent bundle of $M$ \cite{17,\,18}:
$$
\nabla^*\nabla\s\,=\,|\'\nabla\s\'|^2\s,
\tag1.1
$$
where $\nabla$ is the Levi-Civita connection and $\nabla^*\nabla=-\Tr \nabla^2$. The coefficient $|\'\nabla\s\'|^2$, which is (twice) the {\sl vertical energy density\/} of $\s$, is sometimes also called the {\sl bending\/} of $\s$ \cite{17}, and (1.1) are the Euler-Lagrange equations for the associated functional on unit fields.  The study of harmonic unit fields has become popular in recent years, as witnessed by the bibliography of \cite{8}.  There are a number of results specific to the $3$-dimensional case (see \cite{10,\,11,\,12,\,14} amongst others): for example, if $\s$ is a unit Killing vector, or the characteristic vector of a contact metric structure, then $\s$ is harmonic if and only if $\s$ is (pointwise) an eigenvector of the Ricci operator.  The same characterization applies when $\s$ is tangent to the fibres of a {\sl harmonic morphism\/} from $M^3$ to a surface \cite{2}.  Furthermore, in space forms of all dimensions, radial unit vector fields about a point, or a totally geodesically submanifold, are harmonic \cite{3}.  These lead to a number of ``natural'' examples, such as a classification of left-invariant harmonic unit fields on $3$-dimensional unimodular Lie groups, and consequent classification (including stability) of invariant harmonic unit fields on compact quotients.  A unified account of the local and global classification of harmonic morphisms from $3$-dimensional space forms can be found in \cite{2}.  In this paper a number of ``less natural'' examples of harmonic unit fields are constructed, when $M$ is an open subset of $\Bbb R^3$ or $H^3$.   In the Euclidean case we may view $\s$ as a map into $S^2$, and (1.1) is then equivalent to the harmonic map equations \cite{7}; thus, we are constructing harmonic maps $M\to S^2$.  However since this is no longer the case when $M\subset H^3$ we will work with equations (1.1) throughout.  The basic idea is to simplify (1.1) by breaking the symmetry of $M$ and asking $\s$ to respect this additional structure.  The symmetry reduction is guided by the availability of a reasonable supply of simple examples of harmonic unit fields which satisfy the corresponding constraint.  In the Euclidean case we simply ask $\s$ to be {\sl planar;} that is, tangent to a foliation by parallel planes.  We may then (locally) factorise $\s$:
$$
M@>u>>\Bbb R@>w>>S^1\hookrightarrow S^2 
$$
where $w(t)=e^{it}$ and $S^1$ is included in $S^2$ as the equator cut out by the foliation, and this reduces (1.1) to the linear equation $\Delta u=0$.  It is then quite straightforward to restore the full symmetry in order to construct non-constrained harmonic $\s$ (Proposition 1).  One equation in the resultant system of nonlinear PDEs may be viewed as a generalization of the elliptic $3$-dimensional sine-Gordon equation.  In \S3 we use this approach to construct a family of {\it continuous\/} maps $\Bbb R^3\to S^2$, parametrized by an open M\"obius band, all of which are harmonic (hence smooth) on the complement of a line, and whose energy density extends to a bounded continuous function on $\Bbb R^3$.
\par
In \S4 we consider unit vector fields in $H^3$ which are {\sl horospherical;} that is, tangent to a foliation of $M$ by horospheres.  Working with the upper half space model, equation (1.1) then reduces to a linear second order equation (viz\. the hyperbolic Laplace equation) with a nonlinear first order constraint (Proposition 2).  The latter cuts down the availability of solutions quite dramatically compared with the Euclidean case.  Nevertheless, we are able to construct a smooth $1$-parameter family of harmonic unit fields on $H^3$, no member of which is invariant under any discrete subgroup of isometries with compact quotient, but whose energies are arbitrarily close to the minimum for horospherical unit fields, this being achieved only by the invariant unit fields.  In \S5 we examine the stability of unit fields $\s$ whose integral curves are geodesics which are parallel in the hyperbolic sense ie\. convergent to a common point on the ideal boundary of $H^3$.  We refer to such fields as {\sl H-parallel.}  They are harmonic \cite{9}, and their variation fields are precisely the compactly supported horospherical vector fields, for the orthogonal horosphere foliation.  In contrast to the Hopf vector field on the $3$-dimensional sphere \cite{4,\,14}, we show that $\s$ is unstable in the following way: there exists $R_u>0$ such that for any geodesic ball $B_R$ of radius $R>R_u$ it is possible to construct a compactly supported energy-decreasing variation of $\s$ whose support lies inside $B_R$.  We summarise this by saying that $\s$ is {\sl stretchable\/} in all sufficiently large geodesic balls.
\par
It is a pleasure to thank Evgeny Sklyanin, and the referee, for their perceptive and helpful comments.
\bbigskip
\head
2. Generalities, and Euclidean Preliminaries
\endhead
We first record some general identities in Riemannian geometry related to equations (1.1). If $X$ is a vector field on $M$ and $f\colon M\to\Bbb R$ is a smooth function, where $M$ is any Riemannian manifold, then:
$$
\nabla^*\nabla(fX)\,
=\,f\,\nabla^*\nabla X\,-\,2\,\nab{\nabla f}X\,+\,\Delta f\,X,
\tag2.1
$$
where $\nabla f$ is the gradient vector, and $\Delta f$ is the Laplacian:
$$
\Delta f\,=\,-\Tr\nabla df\,
=\,-\divergence(\nabla f).
$$
Note our choice of sign for $\Delta f$.
Furthermore if $g\colon\Bbb R\to\Bbb R$ is smooth then:
$$
\align
\nabla(g(f))\,&=\,g'(f)\,\nabla f,
\tag2.2 \\
\Delta(g(f))\,&=\,g'(f)\,\Delta f\,-\,g''(f)\,|\'\nabla f\'|^2.
\tag2.3
\endalign
$$
\comment
Or if $g\colon M\to\Bbb R$ is smooth then:
$$
\Delta(fg)\,=\,f\'\Delta g\,+\,g\,\Delta f\,-\,2\<\nabla f,\nabla g\>
\tag2.4
$$
\endcomment
\indent
Suppose now that $M$ is a domain of Euclidean $\Bbb R^3$, and $\s$ is a smooth unit vector field on $M$.  Let $(\xi_1,\xi_2,\xi_3)$ be a right-handed orthonormal basis of $\Bbb R^3$, and let $L$ denote the line generated by $\xi_3$, and $P$ the plane spanned by $\xi_1$ and $\xi_2$.  Let $\F$ be the foliation of $M$ by planes parallel to $P$, and let $J$ be the complex structure in $\F$ compatible with the induced orientation.  Thus $\xi_2=J\xi_1$.  If $\s(x)\neq\pm\xi_3$ for all $x\in M$ then $\s$ may be written in {\sl polar form:}
$$
\s\,=\,\cos u\sin v\,\xi_1\,+\,\sin u\sin v\,\xi_2\,+\,\cos v\,\xi_3,
\tag2.4
$$
where $u,v$ are smooth $\Bbb R$-valued functions.  In general, $u$ is only locally defined, any two differing by an integer multiple of $2\pi$.  The {\sl equatorial part\/} $\xi$ of $\s$ is then the unit field:
$$
\xi\,=\,\cos u\,\xi_1\,+\,\sin u\,\xi_2.
$$
Imposing the harmonic section equations (1.1) on $\s$ yields an overdetermined, but consistent, system of nonlinear PDEs for $u$ and $v$.
\proclaim{Proposition 1}
Let $\s$ be a smooth unit vector field on a domain of Euclidean $3$-space, expressible in polar form \rom{(2.4)}.  Then the bending of $\s$ is:
$$
|\'\nabla\s\'|^2\,=\,|\'\nabla v\'|^2\,+\,\sin^2v\,|\nabla u\'|^2,
\tag2.5
$$
and $\s$ is harmonic if and only if:
$$
\align
\Delta u\,&=\,2(\nabla u\dotprod\nabla v)\cot v, 
\tag2.6\\
\Delta v\,&=\,-|\'\nabla u\'|^2\cos v\sin v.
\tag2.7
\endalign
$$
\endproclaim
\demo{Proof}
In terms of its equatorial part $\xi$, we have:
$$
\s\,=\,\sin v\,\xi\,+\,\cos v\,\xi_3.
$$
Now:
$$
\nab X\xi\,=\,-(X.u)\sin u\,\xi_1\,+\,(X.u)\cos u\,\xi_2\,=\,(X.u)J\xi,
$$
and therefore:
$$
\nab X\s\,=\,(X.v)\cos v\,\xi\,+\,(X.u)\sin v\,J\xi\,
-\,(X.v)\sin v\,\xi_3,
$$
from which (2.5) follows.  Furthermore, applying (2.1)--(2.3) yields:
$$
\align
\nabla^*\nabla\xi\,&=\,\Delta(\cos u)\,\xi_1\,+\,\Delta(\sin u)\,\xi_2 \\
&=\,(-\sin u\,\Delta u\,+\,\cos u\,|\'\nabla u\'|^2)\xi_1\,
+\,(\cos u\,\Delta u\,+\,\sin u\,|\'\nabla u\'|^2)\xi_2 \\
&=\,|\'\nabla u\'|^2\xi\,+\,\Delta u\,J\xi,
\endalign
$$
and therefore:
$$
\align
\nabla^*\nabla\s\,
&=\,\sin v\,\nabla^*\nabla\xi\,
-\,2\cos v\,\nab{\nabla v}\xi\,
+\,\Delta(\sin v)\,\xi\,
+\,\Delta(\cos v)\,\xi_3 \\
\vspace{1ex}
&=\,\sin v\,(|\'\nabla u\'|^2\xi\,+\,\Delta u\,J\xi)\,
-\,2(\nabla u\dotprod\nabla v)\cos v\,J\xi \\
&\qquad
+\,(\cos v\,\Delta v\,+\,\sin v\,|\'\nabla v\'|^2)\'\xi\,
+\,(-\sin v\,\Delta v\,+\,\cos v\,|\'\nabla v\'|^2)\'\xi_3 \\
\vspace{1ex}
&=\,|\'\nabla\s\'|^2\s
+\,(\Delta v+|\'\nabla u\'|^2\cos v\sin v)(\cos v\,\xi-\,\sin v\,\xi_3) \\
&\qquad
+\,(\Delta u\'\sin v-2(\nabla u\dotprod\nabla v)\cos v)J\xi,
\endalign
$$
where the final equation is obtained using (2.5).  The result follows on comparison with (1.1), since $(\xi,J\xi,\xi_3)$ is an orthonormal frame.
\qed
\enddemo
\remark{Remarks}
\flushpar
(1)\quad
If $|\'\nabla u\'|^2$ is constant, say $\mu$, then (2.7) is the $3$-dimensional elliptic sine-Gordon equation \cite{5}:
$$
\Delta V\,=\,-\mu\'\sin V,
\qquad\text{where $V=2v$.}
$$
(2)\quad
We say that $\s$ is {\sl planar\/} if it is tangent to $\F$.  In this case $v=\pi/2$, the bending of $\s$ is $|\'\nabla u\'|^2$, and $\s$ is harmonic if and only if $\Delta u=0$.
\flushpar
(3)\quad
The natural $\Bbb C$-action on $\F$, arising from its complex structure, induces a circle-action on planar unit fields.  This extends to a circle-action on unit fields in polar form, when applied to their equatorial part; thus, the action of $e^{it}$ on $\s$ is the unit field $\s_t$ with $u_t=u+t$.  It follows from Proposition 1 that if $\s$ is harmonic then $\s_t$ is a $1$-parameter variation through harmonic unit fields, all of which have identical bending, but which in general are non-congruent. 
\endremark
We illustrate Proposition 1 by first showing how it may be used to reconstruct some known examples of harmonic unit fields \cite{3}.  These will provide useful points of comparison for our forthcoming new examples.
\example{Example 1}
Let $(r,\theta,z)$ be cylindrical polar coordinates with respect to $(\xi_1,\xi_2,\xi_3)$, and let $u=\theta$ and $v=\pi/2$.  Then $\Delta u=0$, and the corresponding harmonic unit field is $\s=\hat\bold r$, the unit radial field about $L$.  The $1$-parameter variation $\s_t$ of $\s$ includes $\s_{\pi/2}=\hat\btheta$, the unit vector in the $\theta$-direction.  All the $\s_t$ are planar, with domain $M=\Bbb R^3\smallsetminus L$, and unbounded bending:
$$
|\'\nabla\s_t|^2\,=\,|\'\nabla\theta\'|^2\,
=\,|\'\hat\btheta/r\'|^2\,=\,\frac{1}{r^2}.
$$
\endexample
\example{Example 2}
Let $(R,\theta,\phi)$ denote spherical polar coordinates relative to
$(\xi_1,\xi_2,\xi_3)$, and let $u=\theta$ and $v=\phi$.
Then $\Delta u=0$, and:
$$
\nabla u\,=\,\hat\btheta/r,\qquad
\nabla v\,=\,\hat\bphi/R.
$$
From the Laplacian in spherical polar coordinates:
$$
\Delta v\,
=\,-\frac{1}{R^2}
\frac{1}{\sin\phi}\frac{\partial}{\partial\phi}
\sin\phi\,
=\,-\frac{1}{R^2}\cot\phi\,
=\,-\frac{\sin\phi\cos\phi}{r^2}.
$$
Therefore both (2.6) and (2.7) hold.  The corresponding harmonic unit field is $\s=\hat\bold R$, the unit radial field about the origin.  The equatorial part $\xi_t$ of $\s_t$ is precisely the $1$-parameter planar family constructed in Example 1; however, $\s_t$ now extends to a harmonic unit field on $\Bbb R^3\smallsetminus\{0\}$.  The $\s_t$ again have unbounded bending:
$$
|\'\nabla\s_t|^2\,=\,\frac{2}{R^2}.
$$
In this case, $\s$ is tangent to the fibres of a harmonic morphism; indeed, when regarded as a map to $S^2$, $\s$ {\it is\/} a harmonic morphism.  The existence of a variation through harmonic maps of identical energy density, none of which is a morphism, is an interesting feature of this example. 
\endexample
\bigskip
\head
3.  New examples in Euclidean space
\endhead
Let $\xi$ denote any member of the $1$-parameter family of planar harmonic unit fields from Example 1, defined on $M=\Bbb R^3\smallsetminus L$.  By a {\sl harmonic desingularization\/} of $\xi$ we mean a harmonic unit field $\s$ on a strictly larger domain $\tilde M$ satisfying:
$$
\s\,-\,(\s\dotprod\xi_3)\xi_3\,=\,\lambda\,\xi,
\tag3.1
$$ 
where $\lambda\colon\tilde M\to\Bbb R$ is a smooth function.  Example 2 may be regarded as a harmonic desingularization of $\xi$ over $\tilde M=\Bbb R^3\smallsetminus\{0\}$.  We now attempt to construct a harmonic desingularization on $\tilde M=\Bbb R^3$.  Continuity of the left hand side of (3.1) implies that $\lambda$  vanishes on $L$.  Now the derivative of $\lambda\,\xi$ is:
$$
d\lambda\otimes\xi\,+\,\lambda\,d\theta\otimes J\xi,
$$
and since the derivative of the left hand side of (3.1) is continuous it follows that $d\lambda$ also vanishes on $L$.  Successive differentiations of $\lambda\,\xi$ lead to the conclusion that all derivatives of $\lambda$ vanish on $L$.  Now, by the regularity of harmonic maps \cite{6}, $\s$ is real analytic when viewed as a function $\Bbb R^3\to S^2$.  Hence the left hand side of (3.1) is also real analytic, viewed as a function $\Bbb R^3\to\Bbb R^2$, and since all its $k$-jets vanish on $L$, it is constant.  Because $\s=\xi_3$ on $L$ it follows that $\s=\xi_3$ identically, and the only harmonic desingularization of $\xi$ is therefore trivial: $\lambda=0$.
\par
We now show how to construct examples of {\it continuous\/} unit fields $\s$ on $\Bbb R^3$ which satisfy (3.1) nontrivially and are harmonic on $M=\Bbb R^3\smallsetminus L$.  Expressing $\s$ in polar form (2.4) on $M$, it follows from (2.6) that $\hat\btheta\dotprod\nabla v=0$. Thus $v$ is rotationally symmetric about $L$, and we may write $v=v(r,z)$ in terms of cylindrical polar coordinates $(r,\theta,z)$.  Suppose for simplicity that $v=v(r)$.  Then:
$$
-\Delta v\,=\,v''\,+\,\frac{v'}{r},
$$
and (2.7) reduces to the following nonlinear ODE with a singularity at $r=0$:
$$
r^2V''+rV'\,=\,\sin V,
\qquad\text{where $V=2v$.}
\tag3.2
$$
In logarithmic coordinates $t=\log r$, (3.2) transforms to:
$$
V''\,=\,\sin V,
$$
and putting $x=\pi-V$ yields the classical pendulum equation:
$$
x''+\sin x\,=\,0.
\tag3.3
$$
Continuity of $\s$ along $L$ requires $v(0+)=0$, which translates to $x(-\infty)=\pi$.  Furthermore if $y=x'$ then:
$$
y(-\infty)\,=\,-\lim_{t\to-\infty}V'(s)\,=\,\lim_{r\to0}rV'(r)\,=\,0.
$$
So in the familiar phase portrait of (3.3) (see for example \cite{1,\,pp\. 90--92}) the relevant trajectories are those on the separatrix which {\it exit\/} the critical point $(\pi,0)$.  These trajectories approach the neighbouring critical points $(-\pi,0)$ and $(3\pi,0)$ as $t\to\infty$, and we deduce that $v=V/2$ increases (resp\. decreases) monotonically to $\pi$ (resp\. $-\pi$) as $r\to\infty$.  Thus as $r\to\infty$ the unit field $\s$ approaches $-\xi_3$, whereas $\s=\xi_3$ on $L$.  The equation of the relevant piece of the separatrix in the phase plane is:
$$
y\,=\,-2\cos(x/2),\quad
-\pi<x<3\pi.
$$
This first order ODE for $x$ has solution:
$$
e^{-t}\,
=\,\cases
c\'(\sec(x/2)+\tan(x/2)),\quad
-\pi<x<\pi, \\
-c\'(\sec(x/2)+\tan(x/2)),\quad
\pi<x<3\pi,
\endcases
$$
where $c>0$ is a constant of integration, and it follows that:
$$
v'(0+)\,=\,-\frac12\lim_{t\to-\infty}e^{-t}\,x'(s)\,
=\,\pm c\lim_{x\to\pi}\big(\sec(x/2)+\tan(x/2)\big)\cos(x/2)\,
=\,\pm2c.
$$
Thus for each $q\in\Bbb R$ there exists a unique smooth solution $v_q$ of
(3.2) on $\Bbb R^+$ with $v_q(0+)=0$ and $v_q'(0+)=q$, noting that if $q=0$
then $v_q=0$ identically.  From (3.3), $v_q''(0+)=0$ for all $q$.  
\par
Now let $\s_{p,q}$ denote the unit field with $u=\theta+p$ and $v=v_q$. 
The translational invariance $\s_{p+2\pi,q}=\s_{p,q}$ is refined by
the glide symmetry $\s_{p+\pi,q}=\s_{p,-q}$, so the $\s_{p,q}$ are in fact
parametrized by an open M\"obius band.  In all cases, $t\mapsto t\'\xi_3$ is
an integral curve of $\s_{p,q}$.  Further qualitative aspects of the flow of
$\s_{p,q}$ may be illustrated by the following two cases.  We refer to the
direction of $\xi_3$ (resp\. $-\xi_3$) as ``up'' (resp\. ``down'').
\smallskip
{\bf Case 1.}
$p=0$.  
The sheaf of planes with axis $L$ is invariant.  Let $\a=|q\'|/2$.  If $q>0$ then the integral curves diverge from $L$ inside the cylinder $C_\a$ with equation $r=\a$, moving upwards until cutting $C_\a$ orthogonally, after which they move downwards and escape.  The overall picture resembles a fountain.  If $q<0$ the flow diagram is obtained by reflection in $P$ and reversing the direction of flow.
\smallskip
{\bf Case 2.}
$p=\pi/2$.  The family of coaxial cylinders $C_r$ ($r>0$) is invariant.  If $q>0$ then inside $C_\a$ the streamlines are right-handed helices which spiral upwards.  The pitch $\rho_r$ of the helices on $C_r$ diverges to infinity as $r\to0$ and converges to zero as $r\to\a-$.  Outside $C_\a$ the integral curves are downward-spiralling left-handed helices, with $\rho_r\to0$ as $r\to\a+$ and $\rho_r\to\infty$ as $r\to\infty$.  On $C_\a$ the integral curves are horizontal anti-clockwise circles (when viewed from ``above'').  If $q<0$ the flow diagram is again obtained by reflection in $P$ and reversing the flow, resulting in a family of upward-spiralling left-handed (resp\. downward-spiralling right-handed) helices inside (resp\. outside) $C_\a$.
\smallskip
Regarding the bending of $\s_{p,q}$, it follows from Proposition 1 that:
$$
|\'\nabla\s_{p,q}\'|^2\,=\,v_q'(r)^2\,+\,\frac{\sin^2v_q(r)}{r^2},
$$
which by a double application of l'H\^opital's rule  converges (to $2q^2$)
as $r\to0$, in contrast to Examples 1 and 2.  Since $v_q'(r)\to0$ as $r\to\infty$,
$|\'\nabla\s_{p,q}\'|^2$ in fact remains bounded over $\Bbb R^3$.  However
the total bending of $\s_{p,q}$ over $\Bbb R^3$ is not finite.  
\par
In all cases, when $\s_{p,q}$ is regarded as a map to $S^2$ its image is
$S^2\smallsetminus\{-\xi_3\}$, in contrast to Example 2 where the
corresponding maps are all surjective.
\par
A more complicated pendulum equation, with damping and variable gravity,
was used by R.T. Smith in \cite{16} to construct ``joins'' of harmonic maps between
spheres.
\bbigskip
\head
4. New examples in hyperbolic space
\endhead
Now let $M$ be a domain in $H^3$.  We denote the hyperbolic metric by $\<\,,\>$, and initially work with the upper half space model, where $\<\,,\>$ is:
$$
ds^2=(dx^2+dy^2+dz^2)/z^2.
$$  
We say that a unit vector field $\s$ is {\sl horospherical\/} if $\s$ is tangent to a foliation $\F$ of $M$ by horospheres.  Without loss of generality, let $\F$ be the standard horosphere foliation, with leaves $z=\const$  Let $(\xi_1,\xi_2,\xi_3)$ be the standard global orthonormal tangent frame:
$$
\xi_1\,=\,z\'\partial_x,\qquad
\xi_2\,=\,z\'\partial_y,\qquad
\xi_3\,=\,z\'\partial_z,
$$
where $\partial_x=\partial/\partial x$ etc.  Then $(\xi_1,\xi_2)$ is a global orthonormal frame of  $\F$, and the integral curves of $\xi_3$ are the unit speed geodesics orthogonal to $\F$.  None of the $\xi_i$ are Killing.  The covariant derivatives of the $\xi_i$ are:
$$
\gathered
\nab{\xi_i}\xi_1\,
=\,\cases 
\xi_3, &i=1, \\
0,&\text{otherwise;} 
\endcases 
\qquad
\nab{\xi_i}\xi_2\,
=\,\cases 
\xi_3, &i=2, \\
0,&\text{otherwise;} 
\endcases \\
\vspace{1ex}
\nab{\xi_i}\xi_3\,
=\,\cases 
-\xi_i, &i=1,2, \\
0,&i=3.
\endcases
\endgathered
\tag4.1
$$
It follows from (4.1) by direct computation that the $\xi_i$ are harmonic, with constant bending:
$$
|\'\nabla\xi_1\'|^2=1=|\'\nabla \xi_2\'|^2,
\qquad
|\'\nabla\xi_3\'|^2=2.
\tag4.2
$$  
The harmonicity of $\xi_i$ was observed in \cite{9}, in the general context of invariant unit fields, where it was also noted that none of the $\xi_i$ are harmonic maps $H^3\to TH^3$.  This may also be seen directly, by inspecting the additional equation \cite{13} for a harmonic unit field $\s$ (on any manifold) to be a harmonic map:
$$
\sum_iR(\s,\nab{E_i}\s)E_i\,=\,0,
$$
where $\{E_i\}$ is any local orthonormal tangent frame.  For a non-Euclidean
space form this yields the characterization that a harmonic unit field $\s$ is a harmonic map if and only if $\s$ is geodesic and solenoidal.  Now $\xi_3$ is geodesic but not solenoidal, whereas $\xi_1$ and $\xi_2$ are solenoidal but not geodesic.  In fact, since the horosphere foliation is not ruled, no horospherical unit field can be a harmonic map.
\par
Any horospherical unit field $\s$ has the following {\sl standard form:}
$$
\s\,=\,\cos u\,\xi_1\,+\,\sin u\,\xi_2,
\tag4.3
$$
where $u$ is a smooth locally-defined $\Bbb R$-valued function.  
\proclaim{Proposition 2}
Let $\s$ be a horospherical unit field, expressed in standard form \rom{(4.3)}.  Then the bending of $\s$ is:
$$
|\'\nabla\s\'|^2\,=\,1+|\'\nabla u\'|^2,
\tag4.4
$$
and $\s$ is harmonic if and only if the following two equations hold:
$$
\align
\Delta u&=0, 
\tag4.5\\
u_x\sin u\,&=\,u_y\cos u,
\tag4.6
\endalign
$$
where $\Delta$ is the hyperbolic Laplacian, and
$u_x=\partial u/\partial x$ etc.  
\endproclaim
\remark{Remark}
The hyperbolic Laplacian in the upper half space model is:
$$
-\Delta u\,=\,z^2(u_{xx}+u_{yy}+u_{zz})\,-\,z\'u_z,
$$
which is of course a linear second order elliptic differential operator.
The nonlinearity of (1.1) appears in the first order constraint equation (4.6).
\endremark
\demo{Proof}
It is convenient to introduce the following almost complex structure
$J$ in $\F$:
$$
J\xi_1\,=\,\xi_2,\qquad
J\xi_2\,=\,-\xi_1.
$$
Then:
$$
\nab X\s\,=\,(X.u)\'J\s\,+\,\cos u\,\nab X\xi_1\,+\,\sin u\,\nab X\xi_2,
$$
and it follows from (4.1) that:
$$
\nab{\xi_1}\s\,=\,(\xi_1.u)\'J\s\,+\,\cos u\,\xi_3,\qquad
\nab{\xi_2}\s\,=\,(\xi_2.u)\'J\s\,+\,\sin u\,\xi_3,\qquad
\nab{\xi_3}\s\,=\,(\xi_3.u)\'J\s.
$$
Therefore:
$$
|\'\nabla\s\'|^2\,
=\,1\,+\,\sum_i\<\nabla u,\xi_i\>^2\,
=\,1\,+\,|\'\nabla u\'|^2.
$$
We now apply (2.1)--(2.3) to (4.3):
$$
\allowdisplaybreaks
\align
\nabla^*\nabla\s\,
&=\,\cos u\,\nabla^*\nabla\xi_1\,
+\,\sin u\,\nabla^*\nabla\xi_2 \\
&\qquad
-\,2\,\nab{\nabla\cos u}\xi_1\,-\,2\,\nab{\nabla\sin u}\xi_2\,
+\,\Delta(\cos u)\,\xi_1\,+\,\Delta(\sin u)\,\xi_2 \\
\vspace{1ex}
&=\,\cos u\,\xi_1\,+\,\sin u\,\xi_2\,+\,2\sin u\,\nab{\nabla u}\xi_1\,
-\,2\cos u\,\nab{\nabla u}\xi_2 \\
&\qquad
-\,\sin u\,\Delta u\,\xi_1\,+\,\cos u\,\Delta u\,\xi_2\,
+\,\cos u\,|\'\nabla u\'|^2\xi_1\,+\,\sin u\,|\'\nabla u\'|^2\xi_2, \\
\intertext{and then apply (4.1) and (4.4):}
&=\,(1+|\'\nabla u\'|^2)\s\,+\,\Delta u\,J\s\,
+\,2(\xi_1.u)\sin u\,\xi_3\,-\,2(\xi_2.u)\cos u\,\xi_3 \\
&=\,|\'\nabla\s\'|^2\s\,+\,\Delta u\,J\s\,
+\,2z(u_x\sin u-u_y\cos u)\'\xi_3.
\endalign
$$
The result follows on comparison with (1.1).
\qed
\enddemo
One consequence of the constraint equation (4.6) is that, in contrast to the Euclidean case, if $\s$ is a harmonic horospherical unit field then the natural circle-action:
$$
e^{it}.\'\s\,=\,\cos(u+t)\xi_1+\sin(u+t)\xi_2
$$
generates a harmonic variation of $\s$ if and only if $u$ is independent of $x$ and $y$; otherwise, the only harmonic members of this $S^1$-family of horospherical unit fields are $\pm\s$.  For example, by analogy with the Euclidean case (Example 1), let $\hat\btheta$ be the following horospherical unit field, defined on $M=H^3\smallsetminus L$ where $L$ is the $z$-axis:
$$
\hat\btheta\,=\,-\frac{y}{\sqrt{x^2+y^2}}\,\xi_1\,
+\,\frac{x}{\sqrt{x^2+y^2}}\,\xi_2.
$$
Then $\hat\btheta$ is a unit Killing field, and therefore harmonic \cite{11}.  However, apart from $-\hat\btheta$, the hyperbolic analogues of the other members of the $1$-parameter family of Example 1 are not harmonic.  In particular, the hyperbolic analogue of $\hat\bold r$ is not harmonic.  This does not contradict the results of \cite{3}, since $\hat\bold r$ is no longer radial. 
\proclaim{Theorem 1}
Let $\s$ be a harmonic horospherical unit field, in standard form \rom{(4.3)}.  
\flushpar
\rom{(i)}\quad
If $u=u(x,y)$ then either $\s$ is invariant, or $\s=\hat\btheta$ upto translation in $H^3$.
\flushpar
\rom{(ii)}\quad
If $u=u(z)$ then $\s$ is a member of the following $2$-parameter family: 
$$
\s_{p,q}\,=\,\cos(pz^2+q)\'\xi_1\,+\,\sin(pz^2+q)\'\xi_2,
\qquad\text{where $p,q\in\Bbb R$.}
$$
\endproclaim
\demo{Proof}
\flushpar
(i)\quad
In this case (4.5) reduces to the Euclidean Laplace equation on $\Bbb
R^2$.  Let $v\colon\Bbb R^2\to\Bbb R$ be the harmonic conjugate of $u$. 
Then
$F=u+iv$ is a holomorphic function of $\zeta=x+iy$.  By the Cauchy-Riemann
equations:
$$
\frac{dF}{d\zeta}\,=\,F_\zeta\,=\,u_x-iu_y,
$$
and hence:
$$
\align
F_\zeta\,e^{iF}\,
&=\,(u_x-iu_y)\,e^{i(u+iv)}\,
=\,e^{-v}(u_x-iu_y)(\cos u+i\sin u) \\
&=\,e^{-v}\big(u_x\cos u+u_y\sin u\,
+\,i(u_x\sin u-u_y\cos u)\big).
\endalign
$$
Therefore (4.6) implies that the holomorphic function $F_\zeta\'e^{iF}$ is purely real, and hence constant.  Consequently we obtain the first order holomorphic ODE:
$$
\frac{dF}{d\zeta}\,=\,k\'e^{-iF},
\qquad
k\in\Bbb R.
$$
This has solution:
$$
F(\zeta)\,=\,-i\log(ik\zeta+\a),
\qquad
\a\in\Bbb C,
$$
with real part:
$$
u(\zeta)\,=\,\arg(ik\zeta+\a).
$$
If $k=0$ then $u$ is constant and $\s$ is therefore invariant.  Otherwise, by
applying a parabolic isometry of $H^3$, we may take $\a=0$, in which case
$\s=\hat\btheta$.
\flushpar
(ii)\quad
In this case equations (4.5) and (4.6) reduce to:
$$
z\'u_{zz}\,=\,u_z,
$$
and it follows that:
$$
u\,=\,pz^2+q,\qquad
p,q\in\Bbb R.
\tag"\qed"
$$
\enddemo
The family $\s_{p,q}$ of Theorem 1 is $2\pi$-periodic in $q$, and hence parametrized by a cylinder; in fact, $q$ is the parameter associated to the natural circle-action. This family of harmonic unit fields is interesting for a number of reasons.  It follows from Proposition 2 that amongst horospherical unit fields those which are invariant minimize energy (over relatively compact domains).  These invariant fields are parallel when restricted to any horosphere, with respect to the induced (Euclidean) metric, and the same is true of the $\s_{p,q}$.  However, by Proposition 2, the bending of $\s_{p,q}$ is:
$$
|\'\nabla\s_{p,q}\'|^2\,=\,1+|\'\nabla u\'|^2\,
=\,1\,+\,z^2(u_x^2+u_y^2+u_z^2)\,
=\,1+4p^2z^4.
$$
It follows that no member of the $1$-parameter family $\s_t=\s_{t,0}$ is invariant under any discrete subgroup of hyperbolic isometries with compact quotient (except when $t=0$), but their local energies are arbitrarily close to the minimum.  At first sight, the existence of the harmonic variation $\s_t$ appears to contradict the variational characterization of harmonic sections.  However on a non-compact manifold harmonic sections are critical points of the vertical energy functional with respect to all {\it compactly supported\/} smooth variations through sections, and the variation field of $\s_t$ is:
$$
-z^2\sin(tz^2)\'\xi_1\,+\,z^2\cos(tz^2)\'\xi_2,
$$ 
which clearly does not have compact support.  
\bigskip
\head
5.  Instability of H-parallel unit fields.
\endhead
In general, a harmonic unit field $\s$ is said to be {\sl unstable\/} if there exists a compactly supported energy-decreasing variation of $\s$.  By \cite{19} the {\sl Jacobi operator\/} for $\s$ is:
$$
\J_\s(\a)\,=\,\nabla^*\nabla\a\,-\,|\'\nabla\s\'|^2\a\,
-\,2\<\nabla\s,\nabla\a\>\s,
$$ 
for all compactly supported smooth vector fields $\a$ pointwise orthogonal
to $\s$.  The {\sl second variation\/} or {\sl Hessian\/} of the
energy/bending functional at $\s$ is then:
$$
\H_\s(\a,\b)\,=\,\int_M\<\J_\s(\a),\b\>\,dV,
$$
where $dV$ denotes the Riemannian volume element, and $\s$ is unstable if
and only if there exists $\a$ with $\H_\s(\a,\a)<0$.  If $U\subset M$ is a
domain containing $\supp(\a)$ then we say that $\s$ is {\sl unstable,} or
{\sl stretchable, in $U$.}  
\proclaim{Theorem 2}
Let $\s$ be an H-parallel unit field on $H^3$.  Then there exists $R_u>0$ such that $\s$ is stretchable in any geodesic ball of radius $R>R_u$.
\endproclaim
\demo{Proof}
Without loss of generality we may assume that $\s=\xi_3$.  Consider the unit field $\xi=\xi_1$, and let $f$ be a smooth function on $H^3$ with compact support.  Then $\a=f\xi$ is a smooth compactly supported variation field for $\s$, and since $|\nabla\s|^2=2$ by (4.2), it follows that:
$$
\<\J_\s(\a),\a\>\,=\,\<\nabla^*\nabla \a,\a\>\,-\,2f^2.
$$
Now $\nabla^*\nabla\xi=\xi$ and it follows from (2.1) that:
$$
\<\nabla^*\nabla\a,\a\>\,=\,\<\nabla^*\nabla(f\xi),f\xi\>\,
=\,\big\<f\,\nabla^*\nabla\xi\,-\,2\,\nab{\nabla f}\xi\,
+\,(\Delta f)\,\xi,\,f\xi\big\>\,
=\,f^2\,+\,f\,\Delta f.
$$
Therefore:
$$
\<\J_\s(\a),\a\>\,=\,f\,\Delta f\,-\,f^2,
$$
and hence by Stokes' Theorem:
$$
\Cal H_\s(\a,\a)\,
=\,\int_M(|\nabla f|^2-f^2)dV.
\tag5.1
$$
\indent
Now for any $x_0\in H^3$ let $B_R(x_0)$ be a geodesic ball about $x_0$ of
radius $R$.  Fix $\d>0$ and suppose that $f$ is a bump function $\v$ with:
$$
\v(x)\,
=\,\cases 
1,& x\in B_R(x_0), \\
0,& x\in M\smallsetminus B_{R+\d}(x_0).
\endcases
$$
Suppose further that $\v$ is radially symmetric about $x_0$, and piecewise linear with $|\nabla\v|=1/\d$ in the spherical shell $S_{R,R+\d}(x_0)=B_{R+\d}(x_0)\smallsetminus B_R(x_0)$.  Such $\v$ may be approximated smoothly, with uniform convergence to both $\v$ and $\nabla\v$.  For further computational simplicity we now work with the unit ball model of $H^3$, where $\<\,,\>$ is:
$$
ds^2\,=\,\frac{4(dx^2+dy^2+dz^2)}{(1-r^2)^2},
$$
with $r^2=x^2+y^2+z^2$.  We also take $x_0=(0,0,0)$.  Then in terms of geodesic spherical polar coordinates $(\rho,\theta,\phi)$ about $x_0$ we have:
$$
\rho\,=\,\tanh(r/2),
$$
and:
$$
\v(x,y,z)\,=\,\v(\rho)\,
=\,\cases
1,& 0\leqs\rho\leqs R, \\
1+\,\dfrac{R-\rho}{\d},&R\leqs\rho\leqs R+1, \\
0,&\rho\geqs R+1.
\endcases
$$
Let $V_R$ (resp\. $V_{R,R+\d}$) denote the volume of $B_R$ (resp\.
$S_{R,R+\d}$), and note that:  
$$
dV\,=\,\sinh^2\rho\,\sin\phi\,d\rho\,d\theta\'d\phi.
$$
It then follows from (5.1) that:
$$
\align
\H_\s(\a,\a)\,
&=\,\frac{1}{\d^2}\int_{S_{R,R+\d}}dV\,
-\,\frac{1}{\d^2}\int_{S_{R,R+\d}}(R+\d-\rho)^2dV\,
-\,\int_{B_R}dV \\
\vspace{1ex}
&=\,\frac{1-(R+\d)^2}{\d^2}V_{R,R+\d}\,-\,V_R\,
+\,\frac{2(R+\d)}{\d^2}\,I_1\,-\,\frac{1}{\d^2}\,I_2,
\endalign
$$
where:
$$
\align
I_1\,&=\,\int_{S_{R,R+\d}}\rho\,dV\,
=\,\int_0^\pi\int_0^{2\pi}\int_R^{R+\d}\rho\sinh^2\rho\sin\phi\,
d\rho\,d\theta\'d\phi \\
\vspace{1ex}
&=\,2\pi\int_R^{R+\d}\rho\'(\cosh2\rho-1)d\rho\,
=\,\pi\left[\rho\sinh2\rho\,-\,\tfrac12\cosh2\rho\,-\rho^2\right]_R^{R+\d}
\endalign
$$
and:
$$
\align
I_2\,&=\,\int_{S_{R,R+\d}}\rho^2dV\,
=\,2\pi\int_R^{R+\d}\rho^2(\cosh2\rho-1)d\rho \\
&=\,\pi\left[(\rho^2+\tfrac12)\sinh2\rho\,
-\,\rho\cosh2\rho\,-\,\frac{2\rho^3}{3}\right]_R^{R+\d}
\endalign
$$
We now note that the volume of a hyperbolic ball of radius $\rho$ is:
$$
V_\rho\,=\,\pi\sinh2\rho\,-\,2\pi\rho,
$$
from which:
$$
\align
I_1\,&=\,(R+\d)V_{R+\d}\,-\,R\'V_R\,
+\,\frac{\pi}{2}\cosh2R\,-\,\frac{\pi}{2}\cosh(2R+2\d)\,
+\,\pi\d(2R+\d), \\
\vspace{1ex}
I_2\,&=\,(R+\d)^2V_{R+\d}\,-\,R^2\'V_R\,
+\,\tfrac12\,V_{R,R+\d}\,
+\,\pi R\cosh2R\,-\,\pi(R+\d)\cosh(2R+2\d) \\
&\qquad
+\,\pi\d\big(4R^2+4R\d+\tfrac43\'\d^2+1\big).
\endalign
$$
Therefore:
$$
\align
2(R+\d)I_1-I_2\,
&=\,(R+\d)^2V_{R+\d}\,-\,R(R+2\d)V_R\,-\,\tfrac12\,V_{R,R+\d} \\
&\qquad
+\,\pi\d\cosh2R\,+\,\pi\d\big(2R\d+\tfrac23\'\d^2-1\big),
\endalign
$$
and so:
$$
\H_\s(\a,\a)\,=\,\left(\frac{1}{\d^2}-\frac12\right)V_{R,R+\d}\,
+\,\frac{\pi}{\d}\cosh2R\,+\,\frac{\pi}{\d}\big(2R\d+\tfrac23\'\d^2-1\big).
$$
Finally:
$$
\align
V_{R,R+\d}\,
&=\,V_{R+\d}-V_R\,
=\,\pi\big(\sinh(2R+2\d)-\sinh2R\big)\,-\,2\pi(R+\d)\,+\,2\pi R \\
&=\,2\pi\cosh(2R+\d)\sinh\d\,-\,2\pi\d,
\endalign
$$
so that:
$$
\H_\s(\a,\a)\,
=\,\frac{\pi}{\d^2}(2-\d^2)\sinh\d\cosh(2R+\d)\,
+\,\frac{\pi}{\d}\cosh2R\,
+\,\frac{\pi}{\d}\big(2R\d+\tfrac53\'\d^2-3\big).
\tag5.2
$$
It follows from (5.2) that if $\d$ is sufficiently large then
$\H_\s(\a,\a)<0$ for all $R>0$.  Let $\d_u$ be the infimum over all
such $\d$, and take $R_u=\d_u$.  \qed
\enddemo
\remark{Remarks}
\flushpar
(1)\quad
From (5.1) we could have estimated:
$$
\H_\s(\a,\a)\,\leqs\,\int_{S_{R,R+\d}}|\'\nabla\v|^2dV\,
-\,\int_{B_R}\v^2dV\,=\,\frac{1}{\d^2}\,V_{R,R+\d}\,-\,V_R.
$$
However, in contrast to Euclidean geometry, the rapid rate of hyperbolic
volume growth renders this upper bound positive for all $R,\d>0$, leaving
us no alternative but to explicitly integrate $\v^2$ over the shell
$S_{R,R+\d}$.
\flushpar
(2)\quad
Further qualitative investigation of (5.2) shows that if $\d$ is 
sufficiently small then $\H_\s(\a,\a)>0$ for all $R>0$.  Let $\d_s$
be the supremum over all such $\d$.  For all $\d_0\in(\d_s,\d_u)$ there
exists $R_0>0$ such that $\H_\s(\a,\a)<0$ for all $R>R_0$.  Numerical
investigation yields estimates $\d_s\approx1.471007$ and $\d_u\approx
1.612195$, and also shows that the transition to instability is
extraordinarily rapid; for example if $\d_0=1.471008$ then
$R_0\leqs8.198206$.  
\flushpar
(3)\quad
It is not known whether the energy functional on unit fields is minimized (over relatively compact domains) by the invariant horospherical fields, or even whether these fields are stable.
\endremark


\Refs
\refstyle{A}
\widestnumber\key{AI}
\ref
\key 1
\by V.\,I. Arnold
\book Ordinary Differential Equations
\publ MIT Press
\yr 1978
\endref
\ref
\key 2
\by P. Baird \& J.\,C. Wood
\book Harmonic Morphisms Between Riemannian Manifolds
\bookinfo  LMS Monographs
\publ OUP
\yr 2003
\endref
\ref
\key 3
\by E. Boeckx \& L. Vanhecke
\paper Harmonic and minimal radial vector fields
\jour  Acta Math\. Hungarica
\yr 2001
\pages 317--331
\endref
\ref
\key 4
\by F. Brito
\paper Total bending of flows with mean curvature correction
\jour Diff. Geom. Appl.
\vol 12
\yr 2000
\pages 157--163
\endref
\ref
\key 5
\by Z. Ding, G. Chen \& S. Li
\paper On positive solutions of the elliptic sine-Gordon equation
\jour  Comm\. Pure \& Applied Analysis
\vol 4
\yr 2005
\pages 283--294
\endref
\ref
\key 6
\by J.\,Eells \& L.\,Lemaire
\paper A report on harmonic maps
\jour Bull\. London Math\. Soc\. 
\vol 10 
\yr 1978
\pages 1--68
\endref
\ref
\key 7
\by J. Eells \& J. H. Sampson
\paper Harmonic mappings of Riemannian manifolds
\jour Amer\. J\. Math\. 
\vol 86
\yr 1964
\page 109--160
\endref
\ref
\key 8
\by O. Gil-Medrano
\paper Unit vector fields that are critical points of the volume and of the
energy: characterization and examples
\inbook Complex, Contact and Symmetric Manifolds
\eds O\. Kowalski, E\. Musso, D\. Perrone
\bookinfo  Progress in Mathematics
\vol 234
\publ Birkhauser
\yr 2005
\pages 165--186
\endref
\ref
\key 9
\by O. Gil-Medrano, J.\,C. Gonz\accute alez-D\accute avila \& L. Vanhecke
\paper Harmonic and minimal invariant unit vector fields on homogeneous
Riemannian manifolds
\jour Houston J\. Math\. 
\vol 27
\yr 2001
\pages 377--409
\endref
\ref
\key 10
\by O. Gil-Medrano \& E. Llinares-Fuster
\paper Second variation of volume and energy of vector fields. Stability of
Hopf vector fields
\jour Math\. Annalen
\vol 320
\yr 2001
\pages 531--545
\endref
\ref
\key 11
\by J.\,C. Gonz\accute alez-D\accute avila \& L. Vanhecke
\paper Minimal and harmonic characteristic vector fields on
three-dimensional contact metric manifolds 
\jour J\. Geom\. 
\vol 72
\yr 2001
\pages 65--76
\endref
\ref
\key 12
\by J.\,C. Gonz\accute alez-D\accute avila \& L. Vanhecke
\paper Energy and volume of unit vector fields on three-dimensional
Riemannian manifolds 
\jour Diff\. Geom\. Appl. 
\vol 16
\yr 2002
\pages 225--244
\endref
\ref
\key 13
\by S-D. Han \& J-W. Yim
\paper Unit vector fields on spheres which are harmonic maps
\jour Math. Z.
\vol 27
\yr 1998
\page 83--92
\endref
\ref
\key 14
\by A. Higuchi, B. S. Kay \& C. M. Wood
\paper The energy of unit vector fields on the 3-sphere
\jour J\. Geom\. Phys\. 
\vol 37
\yr 2001
\page 137--155
\endref
\ref
\key 15
\by D. Perrone
\paper Harmonic characteristic vector fields on contact metric
three-manifolds
\jour Bull\. Austral\. Math\. Soc\.
\vol 67
\yr 2003
\pages 305--315
\endref
\ref
\key 16
\by R.\,T. Smith
\paper Harmonic mappings of spheres
\jour Amer\. J\. Math\. 
\vol 97 
\yr 1975
\page 364--385
\endref
\ref
\key 17
\by G. Wiegmink
\paper Total bending of vector fields on Riemannian manifolds
\jour Math Ann.
\vol 303
\yr 1995
\pages 325--344
\endref
\ref
\key 18
\by C.\,M. Wood
\paper On the energy of a unit vector field
\jour Geom\. Dedicata 
\vol 64
\yr 1997
\pages 319--330
\endref
\ref
\key 19
\by C.\,M. Wood
\paper The energy of Hopf vector fields 
\jour Manuscripta Mathematica
\vol 101 
\yr 2000
\pages 71--88
\endref
\endRefs
\enddocument
\end